\documentclass[12pt]{article} 
\usepackage{latexsym}
\usepackage{amssymb}
\usepackage{euscript}

\let\cal\mathcal

\usepackage{latexsym}
\usepackage{amssymb}
\usepackage{euscript}
\usepackage[dvips]{graphics}
\usepackage{epsf}

\let\cal=\mathcal      

\def\mcc{M\raise.5ex\hbox{c}C}
\def\mccarthy{M\raise.5ex\hbox{c}Carthy}

\def\eg{{\it e.g. }}

\def\l{\lambda}

\def\={\ = \ }

\def\C{\mathbb C}

\def\D{\mathbb D}

\def\NN{\mathbb N}

\def\be{\setcounter{equation}{\value{theorem}} \begin{equation}}
\def\ee{\end{equation} \addtocounter{theorem}{1}}
\def\beq{\begin{eqnarray*}}
\def\eeq{\end{eqnarray*}}
\def\se{\setcounter{equation}{\value{theorem}}} 
\def\att{\addtocounter{theorem}{1}}
\def\vs{\vskip 5pt}
\def\bs{\vskip 12pt}

\def\quest{\bs \att {{\bf Question \thetheorem \ }} }
\def\rem{\bs \att {{\bf Remark \thetheorem \ }} }
\def\bp{{\sc Proof: }}
\def\ep{{}{\hfill $\Box$} \vskip 5pt \par}

\def\bl{\begin{lemma}}
\def\el{\end{lemma}}
\def\bt{\begin{theorem}}
\def\et{\end{theorem}}
\def\bprop{\begin{prop}}
\def\eprop{\end{prop}}
\def\bd{\begin{definition}}
\def\ed{\end{definition}}
\def\br{\begin{remark}}
\def\er{\end{remark}}
\def\bexer{\begin{exercise}}
\def\eexer{\end{exercise}}

\newtheorem{theorem}{Theorem}[section]
\newtheorem{prop}[theorem]{Proposition}
\newtheorem{lemma}[theorem]{Lemma}

\newtheorem{definition}[theorem]{Definition}

\newcommand{\dv}{distinguished variety }
\newcommand{\dvs}{distinguished varieties }
\newcommand{\elt}{\, \in \, }
\newcommand{\umn}{{\cal U}^m_n}
\newcommand{\ut}{U^2_2}

\begin{document}
\setlength{\baselineskip}{21pt}
\title{ Parametrizing Distinguished Varieties}
\author{Jim Agler
\thanks{Partially supported by National Science Foundation Grant
DMS 0400826}\\
U.C. San Diego\\
La Jolla, California 92093
\and
John E. M\raise.5ex\hbox{c}Carthy
\thanks{Partially supported by National Science Foundation Grant
DMS 0501079}\\
Washington University\\
St. Louis, Missouri 63130}

\date{
This paper is dedicated to Joseph Cima on the occasion of his
$70^{\rm th}$ birthday.
}

\bibliographystyle{plain}

\maketitle
\begin{abstract}
A distinguished variety is a variety that exits the bidisk through
the distinguished boundary. We look at the moduli space
for distinguished varieties of rank (2,2).
\end{abstract}

\baselineskip = 18pt

\setcounter{section}{-1}
\section{Introduction}

In this paper, we shall be looking at a special class of bordered
algebraic varieties that are contained in the bidisk $\D^2$ in $\C^2$.
\bd A non-empty 
set $V$ in $\C^2$ is a {\em distinguished variety} 
if there is a polynomial $p$ in $\C[z,w]$ 
such that 
$$
V \= \{ (z,w) \, \in \, \D^2 \ : \ p(z,w) = 0 \} 
$$
and such that
\be
 \overline{V} \cap \partial ( \D^2) \=  \overline{V} \cap (\partial \D)^2
.
\label{eqa1}
\ee
\ed
Condition (\ref{eqa1}) means that the variety exits the bidisk
through the distinguished boundary of the bidisk, the torus.
We shall use $\partial V$ to denote the set given by (\ref{eqa1}):
topologically, it is the boundary of $V$ within the zero set of
$p$, rather than in all of $\C^2$.
\bs
In \cite{agmc_dv}, the authors studied distinguished varieties, 
which we considered interesting because of the following two
theorems:

\bt
\label{thma1}
Let $T_1 $ and $T_2$ be commuting contractive matrices, neither of
which has eigenvalues of modulus $1$.
Then there is a distinguished variety $V$ such that, for any
polynomial $p$ in two variables, the inequality
$$
\| p(T_1, T_2) \| \ \leq \ \| p \|_V 
$$
holds.
\et

\bt
\label{thma2}
The uniqueness variety for a minimal extremal Pick problem on the
bidisk
contains a distinguished variety $V$ that contains each of the
nodes.
\et

It is the goal of this paper to examine the geometry of
distinguished varieties more closely, and in particular
to parametrize the space of all distinguished varieties
of rank $(2,2)$ (see Definition~\ref{defa2} below).

Notice that if $V$ is a distinguished variety, for each $z$ in the
unit disk $\D$, the number of points $w$ satisfying $(z,w) \elt
V$ is constant (except perhaps at a finite number of multiple
points, where the $w$'s must be counted with multiplicity).
So the following definition makes sense:

\bd
\label{defa2}
A \dv is of rank
$(m,n)$ if there are generically $m$ sheets above every
first
coordinate and $n$ above every second coordinate
\ed

The principal result of this paper, Theorem~\ref{thmca},
is a parametrization of 
\dvs of rank $(2,2)$.

\section{Structure theory}
\label{secb}

For positive integers $m$ and $n$, let
\be
\label{eqz05}
U \= \left( \begin{array}{cc} A & B \\
C & D \end{array} \right) \ : \ \C^m \oplus \C^n
\ \to \ \C^m \oplus \C^n
\ee
be an $(m+n)$-by-$(m+n)$ unitary matrix.
Let
\be
\label{eqz1}
\Psi(z) =  A + zB(I- zD)^{-1}C
\ee
be the $m$-by-$m$ matrix valued
function defined on the unit disk $\D$ by the entries of $U$.
This is called the {\it transfer function} of $U$.
Let
$$
U' \= \left( \begin{array}{cc} D^\ast & B^\ast \\
C^\ast & A^\ast \end{array} \right) \ : \ \C^n \oplus \C^m
\ \to \ \C^n \oplus \C^m,
$$
and let
$$
\Psi'(w) \=  D^\ast + wB^\ast(I- w A^\ast)^{-1}C^\ast.
$$
Because $U^\ast U = I$, a calculation yields
\be
\label{eqz2}
I - \Psi(z)^\ast \Psi (z) \=
(1 - |z|^2)\
C^\ast ( I - \bar z D^\ast)^{-1}
(I-z D)^{-1}C ,
\ee
so $\Psi(z)$ is a rational matrix-valued function that is unitary
on the unit circle. Such functions are called rational matrix
inner functions, and it is well-known that all
rational matrix inner functions
have the form (\ref{eqz1}) for some unitary matrix
decomposed as in (\ref{eqz05}) --- see \eg \cite{ampi} for a
proof.
The set
\se\att
\begin{eqnarray}
V &\=&
\{ (z,w) \in \D^2 \ : \
\det(\Psi(z) - w I) = 0 \} 
\label{eqz3}
\se\att
\\
&\=&
\{ (z,w) \in \D^2 \ : \
\det(\Psi'(w) - z I) = 0 \}
\se\att
\\
&\=&
\{ (z,w) \in \D^2 \ : \
\det \left( \begin{array}{cc}
A-wI & zB \\
C & zD - I
\end{array} \right) \= 0 \}
\label{eqz31}
\end{eqnarray}
is a distinguished variety, because when $|z| = 1$, the eigenvalues
of $\Psi(z)$ are unimodular (and a similar statement holds for
$\Psi'$).
The converse was proved in \cite{agmc_dv}: all distinguished
varieties of rank $(m,n)$ can be represented in this way.
So the moduli space for distinguished
varieties of rank $(m,n)$ is a quotient of the space of $(m+n)$-by-$(m+n)$
unitaries. Let us write $\umn$ to denote the set of
$(m+n)$-by-$(m+n)$
unitaries decomposed as in 
(\ref{eqz05}). 
The following result is well-known.

\bprop
\label{propba}
Let $U$ and $U_1$ be in $\umn$, with respective decompositions
$$
U \= \left( \begin{array}{cc} A & B \\
C & D \end{array} \right) \qquad
U_1 \= \left( \begin{array}{cc} A_1 & B_1 \\
C_1 & D_1 \end{array} \right) .
$$
Then they give rise to the same transfer function
iff and only if there is an $n$-by-$n$ unitary $W$ such that
\se\att
\begin{equation}
\label{eqb4}
\left( \begin{array}{cc} A_1 & B_1 \\
C_1 & D_1 \end{array} \right)
\ \= \
\left( \begin{array}{cc} I & 0 \\
0 & W^\ast \end{array} \right)
\
\left( \begin{array}{cc} A & B \\
C & D \end{array} \right)
\
\left( \begin{array}{cc} I & 0 \\
0 & W \end{array} \right)
.
\end{equation}
\eprop
\bp
By looking at the coefficients of powers of $z$ in the transfer
function, we see that $U$ and $U_1$ have the same transfer function
if and only if
\se\att
\begin{eqnarray}
\nonumber
A &\=& A_1 \\
B D^n C &=& B_1 D_1^n C_1 \quad \forall \ n \elt\NN .
\label{eqb5}
\end{eqnarray}
Equation~(\ref{eqb4}) is equivalent to
\beq
A_1 &\=& A \\
B_1 &\=&  BW \\
C_1 &=& W^\ast C \\
D_1 &=& W^\ast D W .
\eeq
Clearly these equations imply (\ref{eqb5}).

To see the converse, note that the fact that $U$ and $U_1$ are
unitaries and $A= A_1$ means $B B^\ast = B_1 B_1^\ast$.
If $B$ is invertible, define $W$
to be $B^{-1} B_1$. This is unitary since $B$ and $B_1$ have the
same absolute values, and then the equations
$BC = B_1 C_1$ and $BDC = B_1 D_1 C_1$ yield
(\ref{eqb4}).

If $B$ is not invertible, then $A$ has norm one. Decompose
$$
\left( \begin{array}{cc} A & B \\
C & D \end{array} \right) \qquad = \qquad
\left(
\begin{array}{cc}
\left(\begin{array}{cc} A' & 0 \\
0 & A'' \end{array} 
\right)
&
\left( \begin{array}{c} 0  \\
B''   \end{array} \right)
\\
\left(\begin{array}{cc} 0 & C''  
\end{array}
\right)
&
D
\end{array} \right),
$$
and apply the same argument to $B''$ and $C''$.
\ep
\rem
Note that $W$ is unique unless $\| A \| = 1$.

\section{Parametrizing distinguished varieties of rank $(2,2)$}
\label{secc}

In this section, we address the question of when two different 
unitaries in $\ut$ give rise to the same \dv. From the previous
section we see that this is equivalent to asking when two
rational matrix inner functions are isospectral.

\bt
\label{thmca}
Let $U, \Psi$ and $V$
be as in (\ref{eqz05}), (\ref{eqz1}) and (\ref{eqz3}),
with $U$ in $\ut$.
Let
$$
U_0 \= \left( \begin{array}{cc} A_0 & B_0 \\
C_0 & D_0 \end{array} \right)
$$
be another unitary in $\ut$. Then $U$ and $U_0$ give rise to the
same \dv iff

\noindent
(i) $A$ and $A_0$ have the same eigenvalues.

\noindent
(ii) $D$ and $D_0$ have the same eigenvalues.

\noindent 
(iii) $BC$ and $B_0 C_0$ have the same trace.
\et

\bp
For simplicity in the proof we will assume that $\det(A) \neq 0$
and that $A$ and $D$ both have two eigenvalues.
(We can attain the remaining cases as a limit of these).

Let
\se\att
\begin{eqnarray}
Q(z,w) &\=& 
\det \left( \begin{array}{cc}
A-wI & zB \\
C & zD - I
\end{array} \right)
\label{eqc1}
\se\att
\\
&\=&
\frac{ \det D}{\det A^\ast}
\
\det \left( \begin{array}{cc}
D^\ast -z & wB^\ast  \\
C^\ast  & wA^\ast - I
\end{array} \right)
\label{eqc2}
\se\att 
\\
&\=& p_2(z) w^2 + p_1 (z) w + p_0 (z) 
\label{eqc3}
\se\att 
\\
&\=& q_2(w) z^2 + q_1 (w) z + q_0 (w) ,
\label{eqc4}
\end{eqnarray}
where $p_{i} $ and $q_j$ are polynomials of degree at most $2$.
As $V$ is the zero set of $Q$, it is sufficient to prove that
conditions (i) --- (iii) completely determine $Q$.
Let $\mu_1$ and $\mu_2$ be the eigenvalues of $D$ and $\l_1$ and
$\l_2$ be the eigenvalues of $A$.

We have 
\beq
p_2(z) &\=& \det (zD - I) \\
&=& (z\mu_1 - 1)(z \mu_2 - 1),
\eeq
so is determined by (ii), the eigenvalues of $D$.
Similarly $q_2(w)$ is determined by (i), the eigenvalues of $A$.

From (\ref{eqc2}) we see that the coefficient of $z^2$ in $Q$
is $( \det D / \det A^\ast)$. Dividing (\ref{eqc3}) by $p_2$, we get
\be
\label{eqc5}
\det ( \Psi(z) - wI) \= w^2 + \frac{p_1(z)}{p_2(z)} w +
\frac{p_0(z)}{p_2(z)} .
\ee
As $\Psi$ is a matrix inner function, we must have that the last
term in (\ref{eqc5}), which is the product of the eigenvalues of
$\Psi$, is inner. Therefore 
$$
p_0(z) \=  e^{i\theta} \
(z - \overline{\mu_1} ) (z - \overline{\mu_2} ).
$$
where 
$$
e^{i\theta} \= (\det D / \det A^\ast) .
$$

It remains to determine $p_1$. 

\bl
\label{lemca}
With notation as above, let 
\be
\label{eqc6}
\det ( \Psi(z) - wI) \= w^2  - a_1(z)  w + a_0(z).
\ee
Then 
\be
\label{eqc7}
a_1(z) \= a_0(z) \overline{ a_1 (\frac{1}{\overline{z}}) } .
\ee
\el
\bp
For any fixed $z$ in $\D$, there are two $w$'s with $(z,w)$ in $V$.
The function $a_0(z)$ is the product of these $w$'s, and $a_1(z)$
is their sum. Labelling them (locally) as $w_1(z)$ and $w_2(z)$,
the right-hand side of (\ref{eqc7}) is
$$
( w_1(z) w_2(z) )\, (\overline{w_1(\frac{1}{\overline{z}}) }
+ \overline{w_2(\frac{1}{\overline{z}}) } ).
$$
When the modulus
of $z$ is $1$, because the variety is distinguished, the right-hand
side of (\ref{eqc7}) equals the left-hand side. By analytic
continuation, they must be equal everywhere.
\ep
Applying the lemma to $- p_1/p_2$ and $p_0/p_2$,
we get
\be
\label{eqc8}
p_1(z) \= 
e^{i \theta} z^2 \overline{p_1(\frac{1}{\overline{z}}) } .
\ee
Writing
$$
p_1(z) \= b_2 z^2 + b_1 z + b_0 ,
$$
(\ref{eqc8}) gives the two equations
\beq
e^{i \theta} \overline{b_2} &\=& b_0 \\
e^{i \theta} \overline{b_1} &\=& b_1 .
\eeq

Comparing (\ref{eqc3}) and (\ref{eqc4}), the coefficient of $z^2 w$
gives us $b_2$ (since we know $q_2$), and hence we also know $b_0$.

Finally, if we know the coefficient of $zw$ in the power series
expansion of (\ref{eqc5}),
we will know 
$$
b_1 p_2(0) - b_0 p_2'(0) ,
$$ and be done.
But
$$
\Psi(z) - wI \= A - wI + zBC + O(z^2) ,
$$
so the coefficient of $zw$ is $-{\rm tr}(BC)$, which is given by
(iii).
\ep

\section{Open Problems}

Two \dvs are 
{\em geometrically equivalent} if 
there is a biholomorphic bijection between
them.

\quest
\label{qa1}
When do two 
unitaries give rise to geometrically equivalent
\dvs?

\vs
Notice that all \dvs of rank $(1,n)$ or $(m,1)$ are 
geometrically equivalent, since they are all biholomorphic
to the unit disk.
\quest
When are two \dvs of rank $(2,2)$ geometrically equivalent?
\vs
\quest
What is 
the generalization of Theorem~\ref{thmca} to \dvs of rank $(2,3)$
or $(3,3)$?
\vs
W.~Rudin showed that smoothly bounded planar domains are
geometrically equivalent to \dvs iff their connectivity is $0$ or
$1$ \cite{rud69b}.
\quest
Which \dvs are geometrically equivalent to planar domains?
\vs
\quest
How can one read the topology of a \dv from a unitary that
determines it as in Section~\ref{secb}?




\bibliography{references}

\end{document}